\documentclass[11pt]{article}
\usepackage{amsmath}
\usepackage{amssymb,amsbsy,amsthm}
\usepackage{graphicx}
\usepackage[dvipsnames]{xcolor}
\usepackage{enumerate}
\usepackage[margin=1in]{geometry}
\usepackage{hyperref}

\allowdisplaybreaks[1]

\numberwithin{equation}{section}

\DeclareMathOperator{\R}{\mathbb{R}} 

\newcommand{\p}{\mathbb{P}} 
\newcommand{\E}{\mathbb{E}} 
\newcommand{\TV}{\mathrm{TV}} 
\newcommand{\Tr}{\mathrm{Tr}} 

\newcommand{\Erf}{\mathrm{Erf}} 

\newtheorem{theorem}{Theorem}[section]

\begin{document}

\title{A smooth transition from Wishart to GOE}
\author{
	Mikl\'os Z.\ R\'acz
	\thanks{Microsoft Research; \texttt{miracz@microsoft.com}.} 
	\and
	Jacob Richey
	\thanks{University of Washington; \texttt{jfrichey@uw.edu}.}
}
\date{\today}

\maketitle


\begin{abstract}
It is well known that an $n \times n$ Wishart matrix with $d$ degrees of freedom is close to the appropriately centered and scaled Gaussian Orthogonal Ensemble (GOE) if $d$ is large enough. Recent work of Bubeck, Ding, Eldan, and Racz, and independently Jiang and Li, shows that the transition happens when $d = \Theta ( n^{3} )$. Here we consider this critical window and explicitly compute the total variation distance between the Wishart and GOE matrices when $d / n^{3} \to c \in (0, \infty)$. This shows, in particular, that the phase transition from Wishart to GOE is smooth. 
\end{abstract}


\section{Introduction} \label{sec:intro} 

The Wishart distribution is a fundamental object appearing in many domains, such as statistics, geometry, quantum physics, and wireless communications, among others. 
In statistics it arises as the distribution of the sample covariance matrix of a sample from a multivariate normal distribution. 
In geometry it is known as the Gram matrix of inner products of $n$ points in $\mathbb{R}^{d}$, 
and it is also the starting point for canonical models of random geometric graphs~\cite{DGLU11,BDER16,EM16}. 

\medskip

It is well known that an $n \times n$ Wishart matrix with $d$ degrees of freedom is close to the appropriately centered and scaled Gaussian Orthogonal Ensemble (GOE) if $d$ is large enough (see, e.g.,~\cite{DGLU11}). 
Recent work~\cite{BDER16,jiang2013approximation} shows that the transition happens when $d = \Theta \left( n^{3} \right)$ and in this paper we study this critical window. 
In Theorem~\ref{thm:smooth} below we explicitly compute the total variation distance between the Wishart and GOE matrices when $d / n^{3} \to c \in (0, \infty)$, 
showing, in particular, that the phase transition from Wishart to GOE is smooth. 

\subsection{Main result} \label{sec:result} 

Let $X$ be an $n \times d$ matrix where the entries are i.i.d.\ standard normal random variables, 
and let $W \equiv W(n,d) = X X^{T}$ be the corresponding $n \times n$ Wishart matrix with $d$ degrees of freedom.\footnote{In statistics the number of samples is usually denoted by $n$ and the number of parameters is usually denoted by $p$, resulting in a $p \times p$ Wishart matrix with $n$ degrees of freedom. Here our notation is taken with the geometric perspective in mind, following~\cite{DGLU11,BDER16,BubeckGanguly15}.} 
Let $M(n)$ be an $n \times n$ matrix drawn from the Gaussian Orthogonal Ensemble, i.e., 
a symmetric $n \times n$ random matrix where the diagonal entries are i.i.d.\ normal random variables with mean zero and variance $2$, 
and the entries above the diagonal are i.i.d.\ standard normal random variables, with the entries on and above the diagonal all independent. 
In order to match the first moment and the scale of the Wishart matrix, we center and scale $M(n)$ appropriately: 
let $M(n,d) := \sqrt{d} M(n) + d I_{n}$, where $I_{n}$ is the $n \times n$ identity matrix. 

\medskip

If $d$ is large enough compared to $n$, then the Wishart matrix becomes approximately like the GOE. 
Recent work of Bubeck, Ding, Eldan, and Racz~\cite{BDER16}, and independently Jiang and Li~\cite{jiang2013approximation}, 
shows that the transition happens when $d = \Theta \left( n^{3} \right)$. 
Specifically, they proved the following theorem, 
where we write $\TV$ for total variation distance.

\begin{theorem}\label{thm:transition}
Define the random matrix ensembles $W(n,d)$ and $M(n,d)$ as above. 
\begin{enumerate}[(a)]
 \item (Bubeck, Ding, Eldan, and Racz~\cite{BDER16}) If $d / n^{3} \to 0$ then 
 \[
  \TV \left( W(n,d), M(n,d) \right) \to 1.
 \]
 \item (Bubeck, Ding, Eldan, and Racz~\cite{BDER16}; Jiang and Li~\cite{jiang2013approximation}) If $d / n^{3} \to \infty$ then 
 \[
  \TV \left( W(n,d), M(n,d) \right) \to 0.
 \]
\end{enumerate}
\end{theorem}

Our focus is on the critical window and our main result is the explicit computation of the limiting total variation distance between $W(n,d)$ and $M(n,d)$ when $d / n^{3} \to c \in \left( 0, \infty \right)$. 

\begin{theorem}\label{thm:smooth}
Define the random matrix ensembles $W(n,d)$ and $M(n,d)$ as above and let $d = d(n)$ be such that $d / n^{3} \to c \in \left( 0, \infty \right)$. 
Then 
\begin{equation}\label{eq:smooth}
\lim_{n \to \infty} \TV \left( W(n,d), M(n,d) \right) = 
\Erf \left( \frac{1}{4\sqrt{3}\sqrt{c}} \right),
\end{equation}
where recall that the error function is defined as 
\[
\Erf \left( x \right) = \frac{2}{\sqrt{\pi}} \int_{0}^{x} e^{-t^{2}} dt. 
\]
\end{theorem}

\begin{figure}[h!]
\centering 
\includegraphics[width=0.4\textwidth]{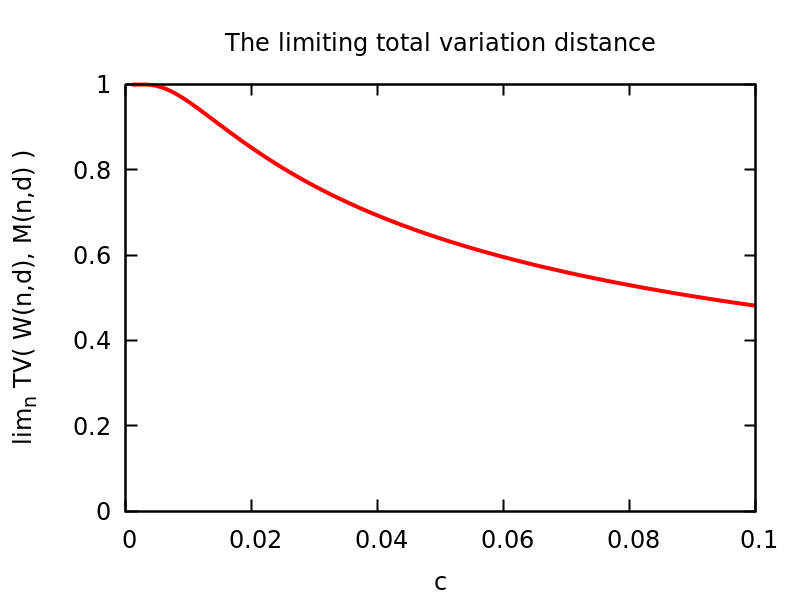}
\caption{The limiting total variation distance as a function of $c$, when $c$ is close to $0$.} 
\label{fig:c_small}
\end{figure}

From this result we can immediately read off that, as $c \to 0$ and $c \to \infty$,  the total variation distance goes to $1$ and $0$, respectively, recovering the previous results described in Theorem~\ref{thm:transition}. 
Since $\Erf \left( x \right) = \tfrac{2}{\sqrt{\pi}} x \left( 1 + o \left( 1 \right) \right)$ as $x \to 0$, the limiting total variation distance decays as 
\[
\Erf \left( \frac{1}{4\sqrt{3}\sqrt{c}} \right) \sim \frac{1}{2 \sqrt{3 \pi} \sqrt{c}}
\]
as $c \to \infty$. 
The behavior of the limit when $c$ is small is plotted in Figure~\ref{fig:c_small}. 

\medskip

From the proof we shall see that the limit in~\eqref{eq:smooth} is the expected value of an explicit function of a two-dimensional Gaussian, 
which comes from the central limit theorem for the first and third moments of the empirical spectral distribution of a GOE matrix.

\subsection{Further related work and open problems} \label{sec:related} 

Several recent works have explored extensions of Theorem~\ref{thm:transition},  
and Theorem~\ref{thm:smooth} raises further questions.

\medskip

\textbf{Robustness.} Bubeck and Ganguly~\cite{BubeckGanguly15} showed that the critical dimension is \emph{universal} in the following sense: Theorem~\ref{thm:transition} holds (up to logarithmic factors) if the entries of $X$ are i.i.d.\ from a sufficiently smooth distribution. 
What can be said about the transition in the critical regime? 
Are there other distributions for which the limiting total variation distance can be computed explicitly? 
If not, can one prove similar qualitative behavior?

\medskip

\textbf{Anisotropy.} Eldan and Mikulincer~\cite{EM16} studied the effect of anisotropy on the power of detecting geometry in random geometric graphs. 
This is directly related to studying Wishart matrices where each row of $X$ is a multivariate normal with a diagonal covariance matrix. 
The authors introduce new notions of dimensionality and prove a theorem similar to Theorem~\ref{thm:transition} with appropriate upper and lower bounds on the ``effective critical dimension''. 
While the primary open problem is to close the gap between these bounds, one may also ask about the nature of the transition at the effective critical dimension: 
can anisotropy cause qualitatively different behavior?

\medskip

\textbf{Other regimes.} 
Theorems~\ref{thm:transition} and~\ref{thm:smooth} state that as $d / n^{3} \to \infty$, \emph{all} statistics of the Wishart $W(n,d)$ and the GOE $M(n,d)$ have asymptotically the same distribution, 
but this is not the case if $d / n^{3}$ remains bounded. 
In the random matrix literature there has been lots of work showing that \emph{particular} statistics of these ensembles have asymptotically the same distribution even when $d \ll n^{3}$. 
For instance, when $d = \Theta \left( n \right)$, then the limiting empirical spectral distribution of the Wishart is the Marchenko-Pastur law, 
which shows the difference between the Wishart and GOE, 
but the largest eigenvalue of the Wishart already behaves like that of the GOE~\cite{johnstone2001distribution,elkaroui2003largest,karoui2007tracy}. 
This naturally raises the question of whether there are other regimes of $d$ and $n$ where there are interesting phase transitions.

\section{Proof of Theorem~\ref{thm:smooth}} \label{sec:proof} 

The main reason that allows for an explicit computation of the limiting total variation distance in Theorem~\ref{thm:smooth} is that both $W \left( n, d \right)$ and $M \left( n, d \right)$ have explicit densities. 
The proof of Theorem~\ref{thm:smooth} is similar to the case of $d / n^{3} \to \infty$ presented in~\cite{BDER16} and proceeds by a Taylor expansion of the ratio of the densities of the two random matrix ensembles. 
The difference compared to the case of $d / n^{3} \to \infty$ is that here the Taylor expansion has to be done to one degree higher. 
As we shall see, taking the limit of the total variation distance as $d / n^{3} \to c \in \left( 0, \infty \right)$ then requires using the central limit theorem for the moments of the empirical spectral distribution of a GOE matrix.

\medskip

\textbf{Step 1: Writing out the total variation distance.} 
Let $\mathcal{P} \subset \R^{\frac{n(n+1)}{2}}$ denote the cone of positive semidefinite matrices. 
It is well known (see, e.g.,~\cite{wishart1928generalised}) that when $d \geq n$, $W(n,d)$ has the following density with respect to the Lebesgue measure on $\mathcal{P}$:
\[
 f_{n,d} \left( A \right) := \frac{\left( \det \left( A \right) \right)^{\frac{1}{2} \left( d - n - 1 \right)} \exp \left( - \frac{1}{2} \Tr \left( A \right) \right)}{2^{\frac{1}{2}dn} \pi^{\frac{1}{4} n \left( n-1\right)} \prod_{i=1}^n \Gamma \left( \frac{1}{2} \left( d+1-i \right) \right)},
\]
where $\Tr\left( A \right)$ denotes the trace of the matrix $A$. 
The density of a GOE random matrix with respect to the Lebesgue measure on $\R^{\frac{n(n+1)}{2}}$ is
$A \mapsto \left( 2 \pi \right)^{-\frac{1}{4} n \left( n + 1 \right)} 2^{-\frac{n}{2}} \exp \left( - \frac{1}{4} \Tr \left( A^2 \right) \right)$ 
and so the density of $M \left( n, d \right)$ with respect to the Lebesgue measure on $\R^{\frac{n(n+1)}{2}}$ is
\[
 g_{n,d} \left( A \right) := \frac{\exp \left( - \frac{1}{4d} \Tr \left( \left( A - d I_n \right)^2 \right) \right)}{\left( 2\pi d \right)^{\frac{1}{4} n \left( n + 1 \right)} 2^{\frac{n}{2}}}.
\]
Denote the measure given by this density by $\mu_{n,d}$, let $\lambda$ denote the Lebesgue measure on $\R^{\frac{n(n+1)}{2}}$ and write $A \succeq 0$ if $A$ is positive semidefinite.
We can then write
\begin{align}
 \TV \left( W(n,d), M \left( n, d \right) \right) &= \int_{\R^{\frac{n(n+1)}{2}}} \left( g_{n,d} \left( A \right) - f_{n,d} \left( A \right) \mathbf{1}_{\left\{ A \succeq 0 \right\}} \right)_{+} d \lambda \left( A \right)  \notag \\
 &= \int_{\R^{\frac{n(n+1)}{2}}} \left( 1 - \frac{f_{n,d} \left( A \right) \mathbf{1}_{\left\{ A \succeq 0 \right\}}}{g_{n,d} \left( A \right)} \right)_{+} d \mu_{n,d} \left( A \right), \label{eq:TV_1_minus_ratio}
\end{align}
where $x_{+} := \max \left\{ x, 0 \right\}$. 
Let $\mathcal{Q}$ denote the set of symmetric matrices for which all of the eigenvalues are in the interval $\left[d - 3 \sqrt{dn}, d + 3 \sqrt{dn} \right]$. 
Since $d / n^{3} \to c > 0$, we have that $\mathcal{Q} \subset \mathcal{P}$ for all $n$ large enough. 
It is known (see, e.g.,~\cite{anderson2010introduction}) that, with probability $1-o\left(1\right)$, 
all the eigenvalues of $M(n)$ are in the interval $\left[-3 \sqrt{n}, 3 \sqrt{n} \right]$, 
which implies that $M\left( n, d \right) \in \mathcal{Q}$. 
Since the integrand in~\eqref{eq:TV_1_minus_ratio} is bounded, we can then write 
\begin{equation}\label{eq:TV}
 \TV \left( W(n,d), M \left( n, d \right) \right) = \int_{\mathcal{Q}} \left( 1 - \frac{f_{n,d} \left( A \right)}{g_{n,d} \left( A \right)} \right)_{+} d \mu_{n,d} \left( A \right)  + o \left( 1 \right)
\end{equation}
and so we may restrict our attention to $\mathcal{Q}$.

\medskip

Define $\alpha_{n,d} \left( A \right) := \log \left( f_{n,d} \left( A \right) / g_{n,d} \left( A \right) \right)$. 
Denote the eigenvalues of an $n \times n$ matrix $A$ by $\lambda_1 \left( A \right) \leq \dots \leq \lambda_n \left( A \right)$; 
when the matrix is obvious from the context, we omit the dependence on $A$. 
Recall that $\det \left(A \right) = \prod_{i=1}^n \lambda_i$ and $\Tr \left( A \right) = \sum_{i=1}^n \lambda_i$. We then have that 
\begin{align*}
 \alpha_{n,d} \left( A \right) &= \frac{1}{2} \sum_{i=1}^n \left\{ \left( d - n - 1 \right) \log \lambda_i - \lambda_i + \frac{1}{2d} \left( \lambda_i - d \right)^2 \right\} \\
 &\quad + \left\{ \frac{n\left( n + 3 \right)}{4} - \frac{dn}{2} \right\} \log 2 + \frac{n}{2} \log \pi + \frac{n\left( n + 1 \right)}{4} \log d  - \sum_{i=1}^n \log \Gamma \left( \frac{1}{2} \left( d + 1 - i \right) \right).
\end{align*}
By Stirling's formula we know that $\log \Gamma \left( z \right) = \left( z - \frac{1}{2} \right) \log z - z + \frac{1}{2} \log \left( 2 \pi \right) + O \left( \frac{1}{z} \right)$ as $z \to \infty$, so 
\begin{align*}
 \alpha_{n,d} \left( A \right) &= \frac{1}{2} \sum_{i=1}^n \left\{ \left( d - n - 1 \right) \log \lambda_i - \lambda_i + \frac{1}{2d} \left( \lambda_i - d \right)^2 \right\} \\
 &\quad + \frac{n\left( n + 1 \right)}{4} \log d  - \frac{1}{2} \sum_{i=1}^n \left( d - i \right) \log \left( d + 1 - i \right) + \frac{1}{2} \sum_{i=1}^n \left( d + 1 - i \right) + O \left( \frac{n}{d} \right).
\end{align*}
Now writing 
$\log \left( d + 1 - i \right) = \log d + \log \left( 1 - \frac{i-1}{d} \right) = \log d  - \frac{i-1}{d}  - \frac{\left( i - 1 \right)^{2}}{2 d^{2}} + O \left( \frac{(i-1)^{3}}{d^3} \right)$ 
we get that
\begin{align*}
 \alpha_{n,d} \left( A \right) = \frac{1}{2} \sum_{i=1}^n \left\{ \left( d - n - 1 \right) \log \lambda_i - \lambda_i + \frac{1}{2d} \left( \lambda_i - d \right)^2 \right\} - \frac{n}{2} \left\{ \left( d - n - 1 \right) \log d - d \right\} - \frac{n^3}{12 d} + o \left( 1 \right). 
\end{align*}
Defining $h \left( x \right) := \frac{1}{2} \left\{ \left( d - n - 1 \right) \log \left( x / d \right) - \left( x - d \right) + \frac{1}{2d} \left( x - d \right)^2 \right\}$, we have that
\begin{equation}\label{eq:alpha_sum_f}
 \alpha_{n,d} \left( A \right)  = \sum_{i=1}^n h\left( \lambda_i \right) - \frac{n^3}{12d} + o \left( 1 \right).
\end{equation}

\medskip

\textbf{Step 2: Taylor expansion and taking the limit.} 
The derivatives of $h$ at $d$ are 
$h \left( d \right) = 0$, 
$h'\left( d \right) = - \frac{n+1}{2d}$, 
$h'' \left( d \right) = \frac{n+1}{2d^2}$, 
$h^{\left( 3 \right)} \left( d \right) = \frac{d - n - 1}{d^3}$, 
$h^{\left( 4 \right)} \left( d \right) = - \frac{3 \left( d - n - 1 \right)}{d^4}$, and also 
$h^{\left( 5 \right)} \left( x \right) = \frac{12 \left( d - n - 1 \right)}{x^{5}}$. 
Approximating $h$ with its fourth order Taylor polynomial around $d$ we get that 
\[
 h(x) = - \frac{n+1}{2d} \left( x - d \right) + \frac{n+1}{4d^2} \left( x - d \right)^2 + \frac{d - n - 1}{6d^3} \left( x - d \right)^3 - \frac{d-n-1}{8 d^4} \left( x - d \right)^4 + \frac{d-n-1}{10 \xi^{5}} \left( x - d \right)^{5},
\]
where $\xi$ is some real number between $x$ and $d$. 
From~\ref{eq:alpha_sum_f} we see that to compute $\alpha_{n,d} \left( A \right)$, 
we need to compute the sum over the eigenvalues $\left\{ \lambda_{i} \right\}_{i=1}^{n}$ of each term in the expansion.

\medskip

First, we argue that the contribution from the remainder term is negligible. 
Recall that $A \in \mathcal{Q}$, and hence $\lambda_{i} \in \left[ d - 3 \sqrt{dn}, d + 3 \sqrt{dn} \right]$ for every $i \in \left[ n \right]$. 
If $x \in \left[ d - 3 \sqrt{dn}, d + 3 \sqrt{dn} \right]$, 
then 
\[
\left| \frac{d-n-1}{10 \xi^{5}} \left( x - d \right)^{5} \right| 
\leq \frac{d-n-1}{10 \left( d - 3 \sqrt{dn} \right)^{5}} \left( 3 \sqrt{dn} \right)^{5} 
= O \left( \frac{1}{n^{2}} \right),
\]
where we used that $d = \Theta \left( n^{3} \right)$. Summing $n$ such terms gives a term of order $O(1/n)$, which is negligible in the limit. 
Turning to the four terms that matter, defining
\begin{equation*}
\begin{aligned}
 S_{1} \left( n, d \right) &:= - \frac{n+1}{2d} \sum_{i=1}^{n} \left( \lambda_{i} - d \right), \\
 S_{3} \left( n, d \right) &:= \frac{d-n-1}{6 d^{3}} \sum_{i=1}^{n} \left( \lambda_{i} - d \right)^{3}, 
\end{aligned}
\qquad \qquad 
\begin{aligned}
 S_{2} \left( n, d \right) &:= \frac{n+1}{4d^{2}} \sum_{i=1}^{n} \left( \lambda_{i} - d \right)^{2}, \\
 S_{4} \left( n, d \right) &:= - \frac{d-n-1}{8 d^{4}} \sum_{i=1}^{n} \left( \lambda_{i} - d \right)^{4}, 
\end{aligned}
\end{equation*}
and also letting 
$S_{0} \left( n, d \right) := - n^{3} / (12 d)$, 
we thus have that 
\begin{equation}\label{eq:alpha_S}
\alpha_{n,d} \left( A \right) = S_{0} \left( n, d \right) + S_{1} \left( n, d \right) + S_{2} \left( n, d \right) + S_{3} \left( n, d \right) + S_{4} \left( n, d \right) + o \left( 1 \right). 
\end{equation}

\medskip

For $i \in \left[ n \right]$ define $\mu_{i} := \tfrac{1}{\sqrt{dn}} \left( \lambda_{i} - d \right)$. 
If $\left\{ \lambda_{i} \right\}_{i=1}^{n}$ are the eigenvalues of $M\left( n, d \right)$, 
then $\left\{ \mu_{i} \right\}_{i=1}^{n}$ are the eigenvalues of $\tfrac{1}{\sqrt{n}} M \left( n \right)$. 
Recall that the empirical spectral distribution $\tfrac{1}{n} \sum_{i=1}^{n} \delta_{\mu_{i}}$ converges weakly to the semicircle distribution with density 
$\rho_{\mathrm{sc}} \left( x \right) = \tfrac{1}{2\pi} \sqrt{4 - x^{2}} \mathbf{1}_{\left\{ \left| x \right| \leq 2 \right\}}$ (see, e.g.,~\cite{anderson2010introduction}). 
With this notation we can rewrite the four quantities above as follows: 
\begin{equation*}
\begin{aligned}
 S_{1} \left( n, d \right) &= - \frac{\sqrt{n} \left( n+1 \right)}{2\sqrt{d}} \times \sum_{i=1}^{n} \mu_{i}, \\
 S_{3} \left( n, d \right) &= \frac{d-n-1}{6 d} \times \sqrt{\frac{n^{3}}{d}} \times \sum_{i=1}^{n} \mu_{i}^{3},
\end{aligned}
\qquad \qquad 
\begin{aligned}
 S_{2} \left( n, d \right) &= \frac{n^{2} \left( n+1 \right)}{4d} \times \frac{1}{n} \sum_{i=1}^{n} \mu_{i}^{2}, \\
 S_{4} \left( n, d \right) &= - \frac{d-n-1}{8d} \times \frac{n^{3}}{d} \times \frac{1}{n} \sum_{i=1}^{n} \mu_{i}^{4}.
\end{aligned}
\end{equation*}

\medskip

Let $U$ be a random variable distributed according to the semicircle law. 
We know that for any fixed $k \in \mathbb{N}$, the $k^{\text{th}}$ moment of the empirical spectral distribution converges in probability to the $k^{\text{th}}$ moment of the semicircle law (see~\cite[Lemmas 2.1.6 and 2.1.7]{anderson2010introduction}), i.e., 
\[
 \frac{1}{n} \sum_{i=1}^{n} \mu_{i}^{k} \to \E \left[ U^{k} \right].
\]
Since $\E \left[ U^{2} \right] = 1$ and $\E \left[ U^{4} \right] = 2$, we have that, as $d / n^{3} \to c \in \left( 0, \infty \right)$, 
$S_{2} \left( n, d \right) \to 1 / (4c)$ and 
$S_{4} \left( n, d \right) \to - 1 / (4c)$, 
so put together we have that 
$S_{2} \left( n, d \right) + S_{4} \left( n, d \right) \to 0$.

\medskip

The central limit theorem for the moments of the empirical spectral distribution of a GOE matrix (see~\cite[Theorem~2.1.31 and Exercise~2.1.35]{anderson2010introduction} and~\cite{anderson2006clt}) shows that  
\[
\left( \sum_{i=1}^{n} \mu_{i}, \sum_{i=1}^{n} \mu_{i}^{3} \right) \Longrightarrow \left( N_{1}, N_{3} \right), 
\]
where $\left( N_{1}, N_{3} \right)$ are jointly normal 
with mean zero and covariance matrix 
$C = 
\left(
\begin{smallmatrix} 
2 & 6 \\ 
6 & 24 
\end{smallmatrix}
\right)$. 
The entries of the covariance matrix $C$ are special cases of the more general formulas found in~\cite{anderson2010introduction,anderson2006clt}, 
so we do not describe the computations here, but one can verify these numbers by using the identity
\[
 \sum_{i=1}^{n} \mu_{i}^{k} = \Tr \left( \left( \tfrac{1}{\sqrt{n}} M \left( n \right) \right)^{k} \right) = \sum_{i=1}^{n} \left( \left( \tfrac{1}{\sqrt{n}} M \left( n \right) \right)^{k} \right)_{i,i}
\]
and computing appropriate moments of normal random variables.

\medskip

Putting everything together we see that, as $d / n^{3} \to c \in \left( 0, \infty \right)$, we have that 
\begin{equation}\label{eq:weak_convergence_moments}
\left( S_{0}, S_{1}, S_{2}, S_{3}, S_{4} \right) \Longrightarrow \left( - \frac{1}{12c}, - \frac{1}{2 \sqrt{c}} N_{1}, \frac{1}{4c}, \frac{1}{6 \sqrt{c}} N_{3}, - \frac{1}{4c} \right). 
\end{equation}
Therefore, since the function 
$\left( s_{0}, s_{1}, s_{2}, s_{3}, s_{4} \right) \mapsto \left( 1 - \exp \left( s_{0} + s_{1} + s_{2} + s_{3} + s_{4} \right) \right)_{+}$ 
is continuous and bounded, we have by~\eqref{eq:TV},~\eqref{eq:alpha_S}, and~\eqref{eq:weak_convergence_moments} that 
\begin{equation}\label{eq:limit_expectation}
 \TV \left( W(n,d), M \left( n, d \right) \right) \to \E \left[ \left( 1 - \exp \left( - \frac{1}{12 c} - \frac{1}{2\sqrt{c}} N_{1} + \frac{1}{6 \sqrt{c}} N_{3} \right) \right)_{+} \right].
\end{equation}

\medskip

\textbf{Step 3: Evaluating the limit.} 
What remains is to evaluate the expectation on the right hand side of~\eqref{eq:limit_expectation}. 
Let $Y$ and $Z$ be independent normal random variables with mean zero and variances $2$ and $6$, respectively. 
Then $\left( N_{1}, N_{3} \right)$ and $\left( Y, 3Y + Z \right)$ have the same distribution, since both are Gaussian with the same mean vector and covariance matrix. 
Notice that 
\[
-\frac{1}{2\sqrt{c}} Y + \frac{1}{6 \sqrt{c}} \left( 3 Y + Z \right) = \frac{1}{6 \sqrt{c}} Z 
\]
and hence the right hand side of~\eqref{eq:limit_expectation} is equal to 
\[
\E \left[ \left( 1 - \exp \left( - \frac{1}{12c} + \frac{1}{6\sqrt{c}} Z \right) \right)_{+} \right].
\]
Since 
$1 - \exp \left( - 1 / (12c) + z / (6\sqrt{c}) \right) \geq 0$ 
if and only if 
$z \leq 1 / (2\sqrt{c})$, 
we have that 
\begin{align*}
\E \left[ \left( 1 - \exp \left( - \frac{1}{12c} + \frac{1}{6\sqrt{c}} Z \right) \right)_{+} \right]
&= \int_{-\infty}^{\frac{1}{2\sqrt{c}}} \left( 1 - e^{- \frac{1}{12c} + \frac{1}{6\sqrt{c}} z} \right) \cdot \frac{1}{\sqrt{12\pi}} e^{-\frac{z^{2}}{12}} dz \\
&= \frac{1}{\sqrt{12\pi}} \int_{-\infty}^{\frac{1}{2\sqrt{c}}} e^{-\frac{z^{2}}{12}} dz - \frac{1}{\sqrt{12\pi}} \int_{-\infty}^{\frac{1}{2\sqrt{c}}} e^{-\frac{\left(z - 1/\sqrt{c} \right)^{2}}{12}} dz \\
&= \p \left( Z < \frac{1}{2\sqrt{c}} \right) - \p \left( Z < - \frac{1}{2\sqrt{c}} \right) = \Erf \left( \frac{1}{4\sqrt{3}\sqrt{c}} \right),
\end{align*}
which concludes the proof.


%


\bibliographystyle{abbrv}
\bibliography{bib}




\end{document}